\documentclass[12pt]{amsart}

\usepackage[english]{babel}
\usepackage[alphabetic]{amsrefs}
\usepackage{graphicx}
\graphicspath{ {images/} }
\usepackage[T1]{fontenc}
\usepackage{txfonts}
\usepackage{times}
\usepackage{amssymb,verbatim,amscd,amsfonts,amsmath,amsthm,enumerate}
\usepackage[all]{xy}
\usepackage{subfig}
\usepackage{tikz}
\usetikzlibrary{shapes,arrows,shadows}
\usetikzlibrary{decorations.markings}
\usepackage{color}
\usepackage[normalem]{ulem}
\usepackage{hyperref}
\usepackage{wrapfig}
\usetikzlibrary{arrows}
\usepackage{pb-diagram}

\newtheorem{theorem}{Theorem}[section]

\newtheorem{remark}[theorem]{Remark}

\begin{document}




\title{On the Infinite Loch Ness monster}


\author{John A. Arredondo}
\email{alexander.arredondo@konradlorenz.edu.co}
\address{Fundaci\'on Universitaria Konrad Lorenz.
CP. 110231, Bogot\'a, Colombia.}

\author{Camilo Ram\'irez Maluendas}
\email{camilo.ramirezm@konradlorenz.edu.co}
\address{Fundaci\'on Universitaria Konrad Lorenz.
CP. 110231, Bogot\'a, Colombia.}

\maketitle

\begin{abstract}
In this paper we present in a topological way the construction of the orientable surface with only one end and infinite genus, called \emph{The Infinite Loch Ness Monster}. In fact, we introduce a flat and hyperbolic construction of this surface. We discuss how the name of this surface has evolved and how it  has been historically understood.
\end{abstract}



\section{Introduction}
\label{Introduction}

\begin{figure}[h!] \label{real-mons}
\begin{center}
  \includegraphics[scale=0.4]{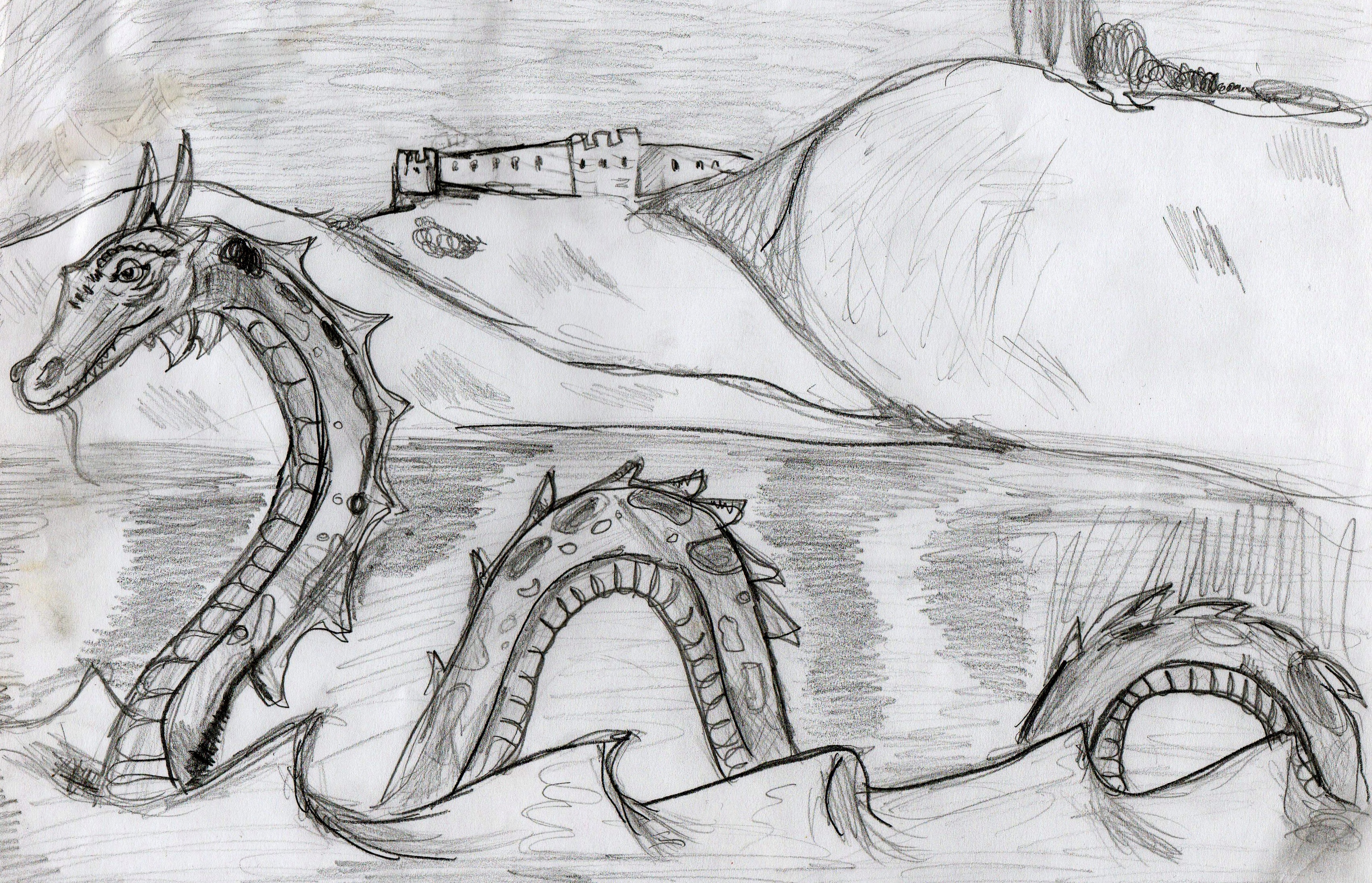}
    \caption{\textit{ Loch Ness Monster in The Great Glen in the Scottish.}}
    \emph{Image by xKirinARTZx, taken from devianart.com}
\end{center}
\end{figure}

The term Loch Ness Monster is  well known around the world, specially in The Great Glen in the Scottish highlands, a rift valley which contains three important lochs for the region, called Lochy, Oich and Ness. The last one, people believe that  a  monster lives and lurks, baptized with the name of the loch. The existence of the monster is not farfetched, people say, taking into account that the Loch Ness  is deeper than the North Sea and is very long, very narrow and has never been known to freeze (see Figure \ref{real-mons}).

The earliest report of such a monster appeared in the Fifth century, and from that time different versions about the monster  passed from generation to generation \cite{Ste}. A kind of modern interest in the monster was sparked by 1933 when George Spicer and his wife saw the monster crossing the road in front of their car. After that sighting, hundreds of different reports about the monster have been collected, including photos, portrayals and other descriptions.  In spite of this evidence, without a body, a fossil or the monster in person, The  Loch Ness Monster is only part of the  folklore.

In a different context, in mathematics,  the term \emph{Loch Ness Monster} is  well known, and not in the folklore. In number theory there is a family of functions called exponential sums, which in general take the form
\begin{equation}
s(n)=\sum_{n=1}^N e^{2\pi i f(n)},
\end{equation}
and for the special case in which
\begin{equation}
f(n)=(ln (n))^4
\end{equation}
the graph of the curve associated to that  function is called \emph{Loch Ness Monster}, dubbed to the curve by J. H. Loxton \cite{Lox}, \cite{Lox1}.

\begin{figure}[h!] \label{curve-mons}
\begin{center}
 \includegraphics[scale=0.5]{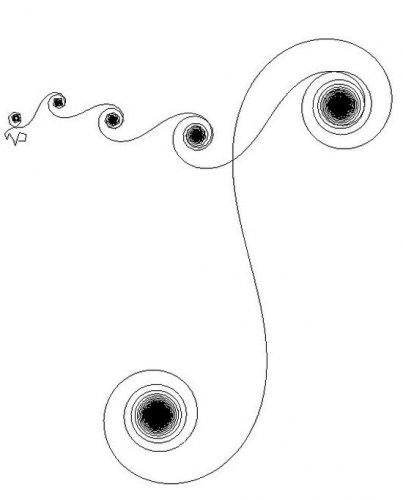}
 \caption{\textit{ Loch Ness Monster curve depicted with $N= 6000$.}}
\end{center}
 \end{figure}

\begin{figure}[h!]\vspace{1cm}
\begin{center}
\includegraphics[scale=0.5]{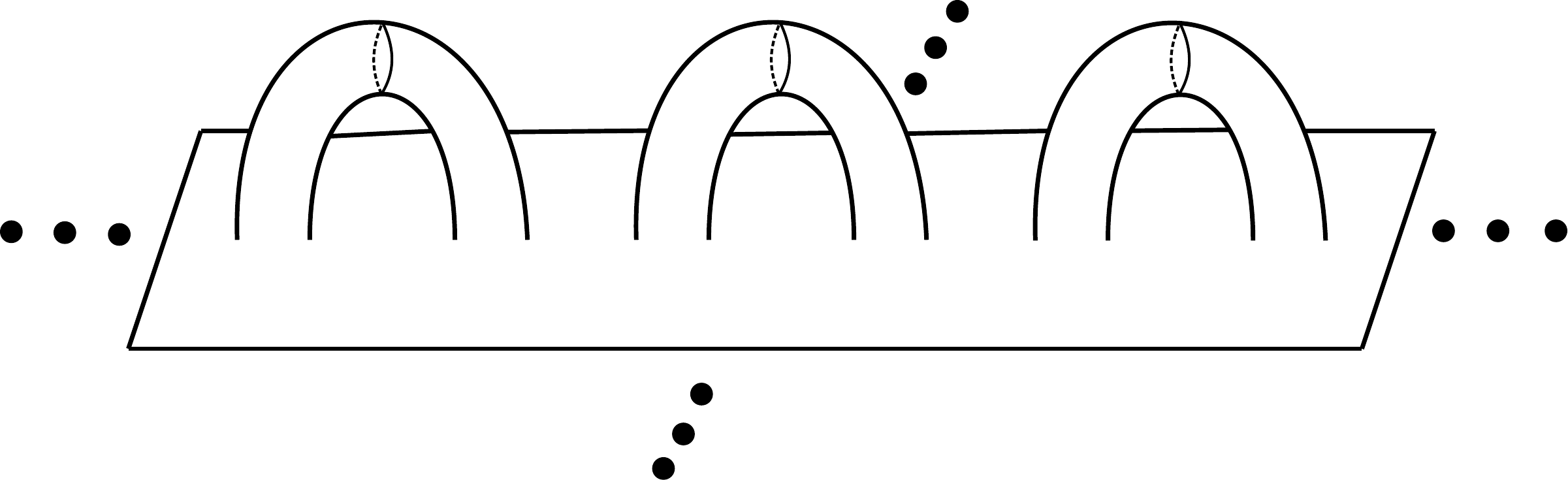}\\
  \caption{\emph{The Infinite Loch Ness monster.}}
   \label{Figure3}
\end{center}
\end{figure}

From view of the Ker\'ekj\'art\'o's theorem of classification of noncompact surfaces (\emph{e.g.}, \cite{Ker}, \cite{Ian}), the \emph{ Infinite Loch Ness Monster} is the name of the orientable surface which has infinite genus and only one end \cite{Val}. Simply, \'E. Ghys (see \cite{Ghy}) describes it as the orientable surface obtained from the Euclidean plane which is attached to an infinity of handles (see Figure \ref{Figure3}). Or alternatively, from a geometric viewpoint one can think that the Infinite Loch Ness monster is the only orientable surface having infinitely many  handles and only one way to go to infinity.

\begin{figure}[h!]
     \centering
     \begin{tabular}{ccc}
      \includegraphics[scale=0.5]{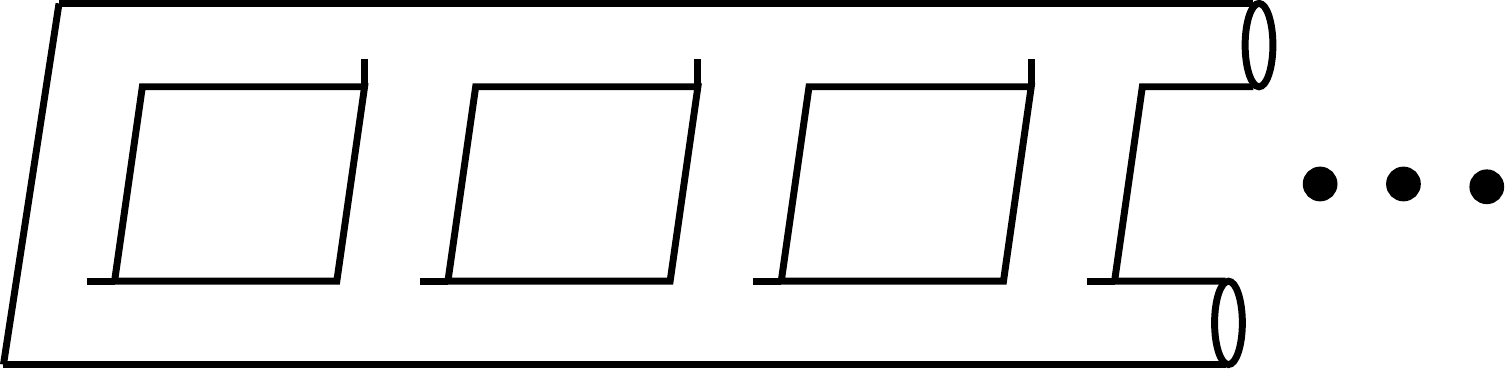} &    & \includegraphics[scale=0.35]{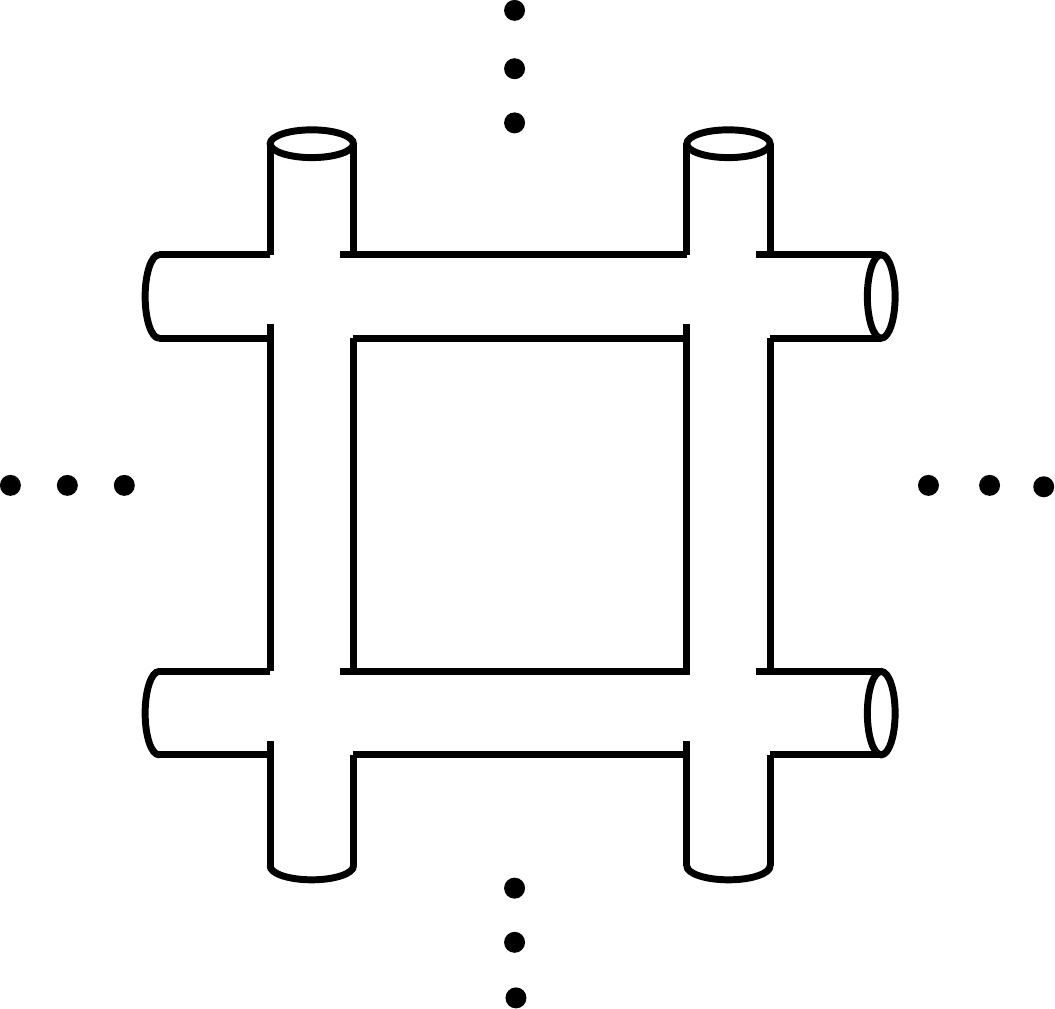} \\
     \text{ \emph{a. Jacob's ladder.}}          &    & \text{\emph{b. Infinite jail cell windows.}} \\
     \includegraphics[scale=0.5]{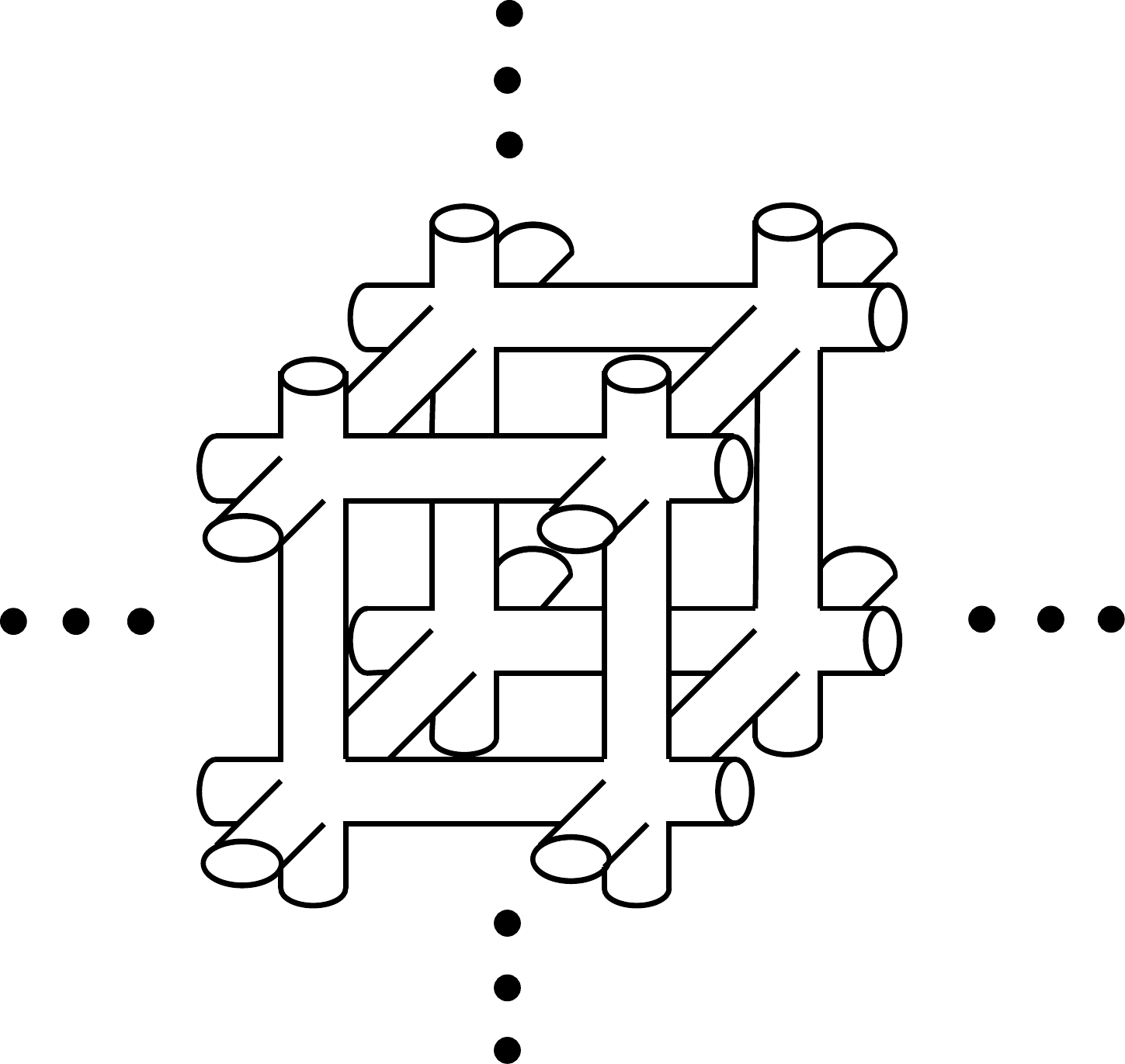} &    & \\
     \text{\emph{c. Infinite jangle gym.}}       &    &
     \end{tabular}
     \caption{\emph{Surfaces having only one end and infinite genus.}}
     \label{Figure4}
\end{figure}

In the seventies, the interest by several authors (\emph{e.g.}, \cite{Sow}, \cite{Ni}, \cite{Can}) on the qualitative study in the noncompact leaves in foliations of closed manifolds had grown. Ongoing in this line of research considering closed 3-manifolds foliated by surfaces, A. Phillips and D. Sullivan proved that the quasi-isometry types of the surfaces well known as the \emph{Jacob's ladder}\footnote{E. Ghys calls Jacob's ladder to the surface with two ends and each ends having infinite genus (see \cite{Ghy}). However, M. Spivak calls this surface the doubly infinite-holed torus (see \cite[p.24]{Spiv})}, the \emph{Infinite jail cell windows} \cite[p.24]{Spiv}, and the \emph{Infinite jangle gym} (see Figure \ref{Figure4}) cannot occur in foliations of $S^3$, or in fact in orientable foliation of any manifold with second Betti number zero. Nevertheless, all these surfaces are diffeomorphic to the \emph{Infinite Loch Ness monster} (see \cite{PSul}). Roughly speaking from the historical point of view, this nomenclature to this topological surface appeared published on \emph{Leaves with Isolated ends in Foliated 3-Manifolds} (\cite[1977]{Can2}), however the authors wedge the term \emph{Infinite Loch Ness monster} to preliminary manuscript of \cite{PSul}, which was published the following year. Under these evident, one can consider to A. Phillips and D. Sullivan as the \emph{Infinite Loch Ness monster}'s parents.

\begin{figure}[h!]
\begin{center}
\includegraphics[scale=0.9]{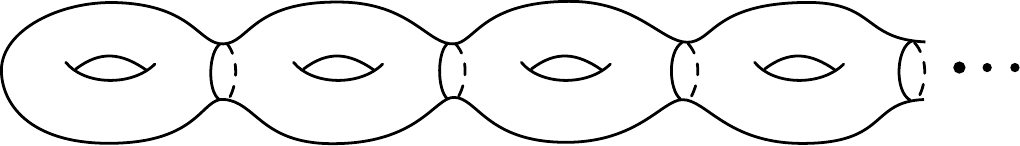}\\
  \caption{\emph{The infinite-holed torus.}}
   \label{Figure5}
\end{center}
\end{figure}

\begin{remark}
Perhaps the reader has found on the literature other names for this surface  with infinite genus and only one end, for example, the \emph{infinite-holed torus}  (see \cite[p.23]{Spiv}). Figure \ref{Figure5}.
\end{remark}

\begin{figure}[h!]
\begin{center}
\begin{tabular}{ccccc}
\includegraphics[scale=0.5]{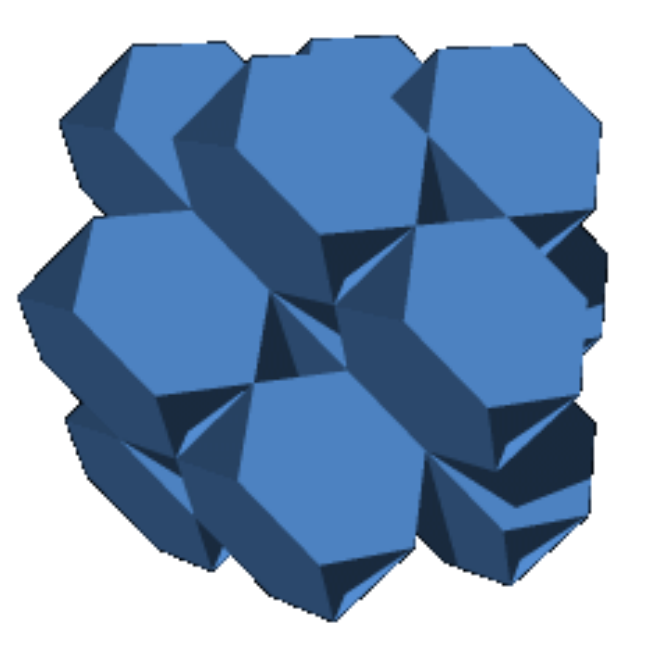}   & \includegraphics[scale=0.26]{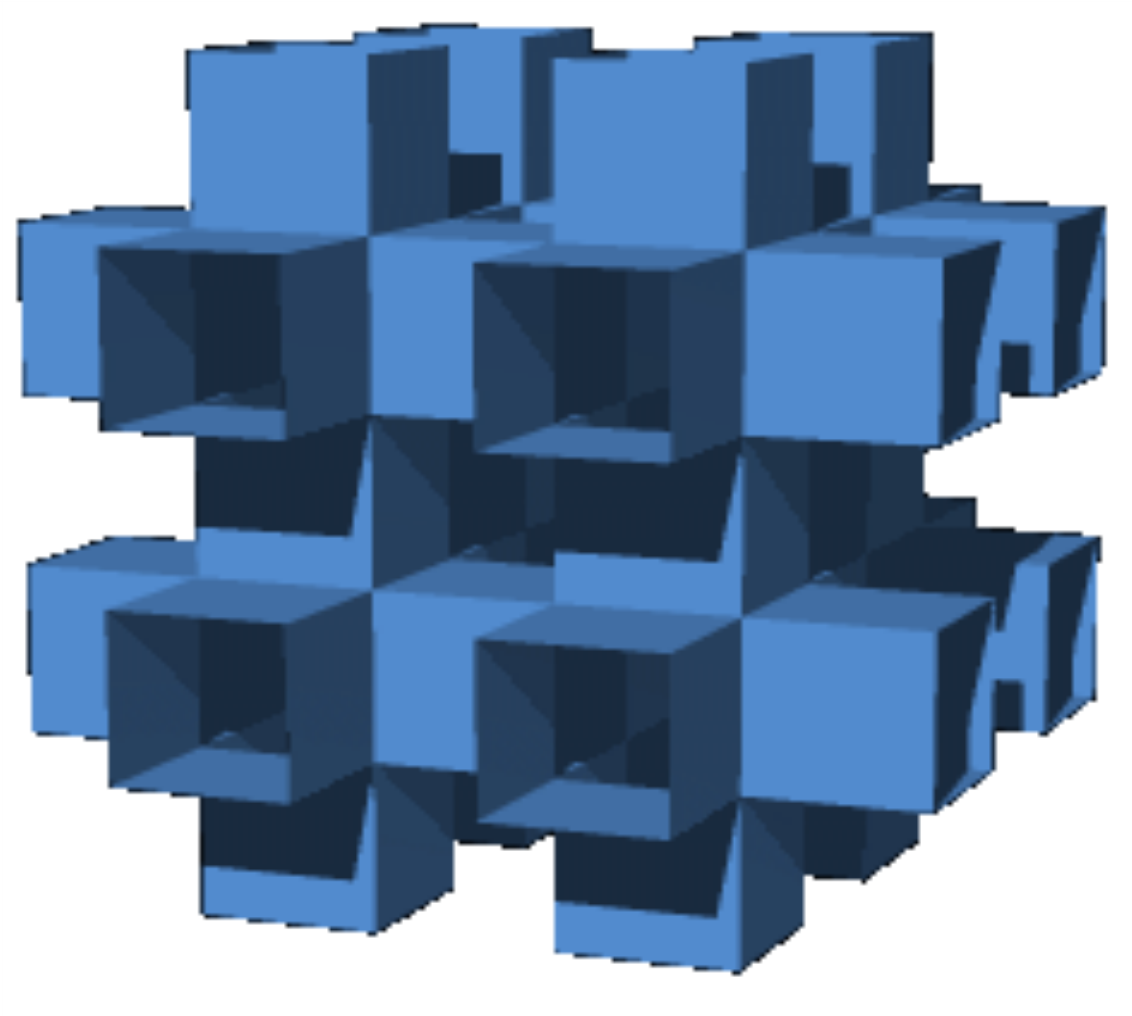}       \\
\text{\emph{a. The multiplied tetrahedron $\mu T$.}}     & \text{\emph{b. The multiplied cube $\mu C$.}}           \\
\includegraphics[scale=0.5]{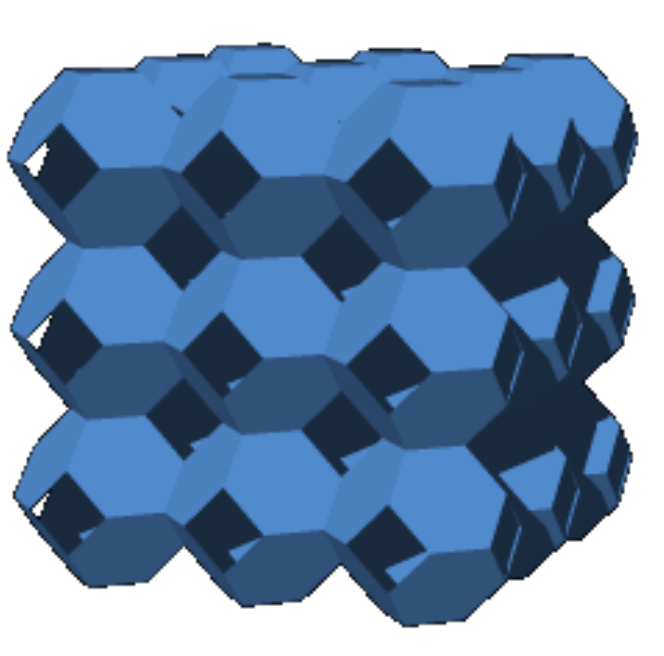}    &\\
\text{\emph{c. The multiplied octahedron $\mu O$.}}      &
\end{tabular}
  \caption{\emph{Locally the skew polyhedra or Coxeter-Petrie polyhedra.}}
  \emph{ Images by Tom Ruen, distributed under CC BY-SA 4.0.}
  \label{Polyhedra}
\end{center}
\end{figure}

The Infinite Loch Ness monster has also appeared in the area of Combinatory. Its arrival was in 1929 when J. P. Petrie told H. S. M. Coxeter that had found two new infinite regular polyhedra. As soon as J. P. Petrie begun to describe them and H. S. M. Coxeter understood this, the second pointed out a third possibility. Later they wrote a paper calling this mathematical objets the \emph{skew polyhedra} \cite{Cox1}, or also known today as the Coxeter-Petrie polyhedra. Indeed, they are topologically equivalent to the Infinite Loch Ness monster as shown in \cite{ARV}.  Given that from a combinatory view one can think that skew polyhedra are multiple covers of the first three Platonic solids, J. H. Conway and \emph{et. al.} \cite[p.333]{Con} called them the \emph{multiplied tetrahedron}, the \emph{multiplied cube}, and the \emph{multiplied octahedron}, and denoted them $\mu T$, $\mu C$, and $\mu O$, respectively. See Figure \ref{Polyhedra}.

\begin{figure}[h!]
\begin{center}
\includegraphics[scale=0.7]{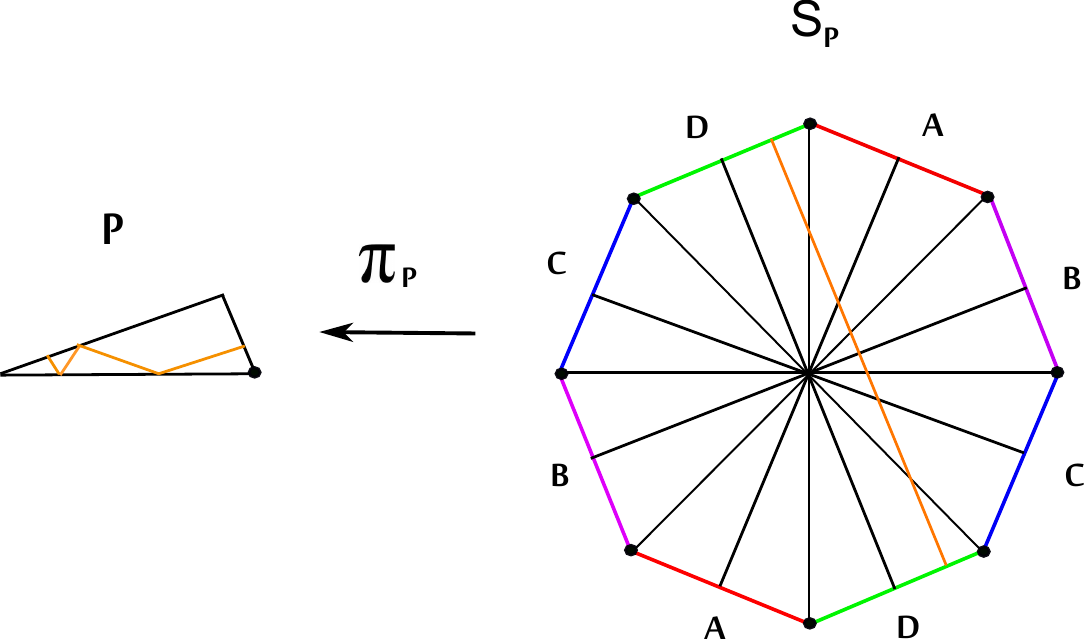}\\
\caption{\emph{Billiard associated to a rectangle triangle.}}
 \label{billar}
\end{center}
\end{figure}

In billiards, an interesting area of Dynamical Systems, during 1936 the mathematicians
\emph{R. H. Fox} and \emph{R. B. Kershner} \cite{Fox} (later, used it by A. B. Katok and A. N. Zemljakov \cite{KZ}) associated to each \emph{billiard} $\phi_P$ coming from an Euclidian compact polygon $P\subseteq \mathbb{E}^2$ a surface $S_{P}$ with structure of translation, which they called \emph{Ueberlagerungsfl\"ache} and means \emph{covered surface}, and a projection map $\pi_p: S_p\to \phi_P$ mapping each geodesic of $S_P$ onto a \emph{billiard trajectory} of $\phi_P$ (see Table \ref{tabla1} and Figure \ref{billar}). Later, F. Valdez published a paper \cite{Val}, which proved that the surface \emph{Ueberlagerungsfl\"ache} $S_P$ associated to the billiard $\phi_P$, being $P\subseteq \mathbb{E}^2$ a polygon with almost an interior vertex of the form $\lambda \pi$ such that $\lambda$ is a irrational number, is the Infinite Loch Ness monster.

\begin{table}[h]
\[
\xymatrix{
                                    & *+[F]{ P\subseteq \mathbb{E}^2 }\ar@(u,u)[dr]^{\text{ \,\, \emph{Ueberlagerungsfl\"ache}}} \ar@(u,u)[dl]_{Billiard}&                                                                   \\
\phi_P   &                                          & S_P\ar[ll]_{\pi_P}
}
\]
\caption{\emph{Billiard $\phi_P $ and surface $S_p$ associated to the polygon $P$.}}
 \label{tabla1}
\end{table}

\section{Building the Infinite Loch Ness Monster}
\label{Building}

\subsection{A tame Infinite Loch Ness Monster}
\label{tame}

\begin{figure}[ht!]
\begin{center}
   \includegraphics[scale=0.6]{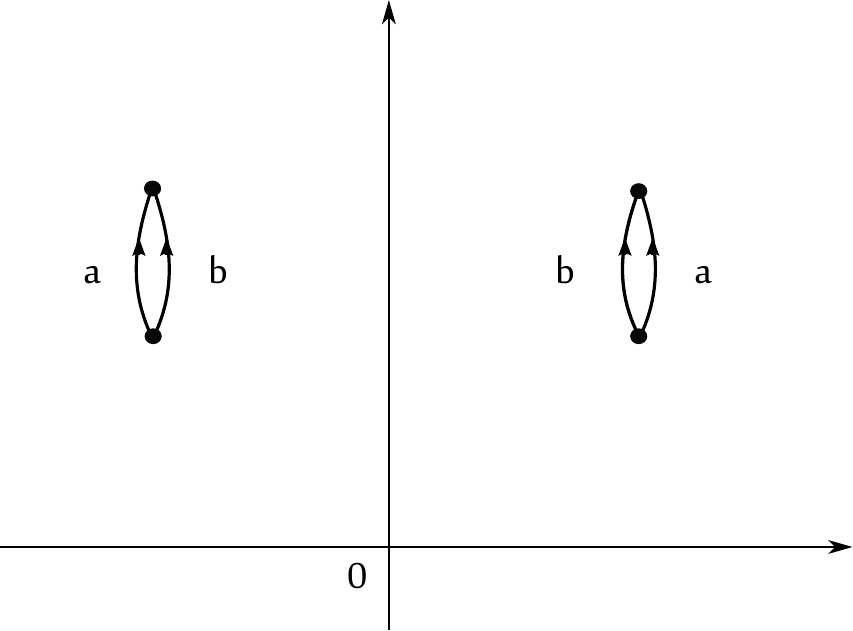}
   \caption{\emph{Two straight segments on $\mathbb{E}^2$.}}
     \label{gluemarks}
 \end{center}
 \end{figure}

An easy and simple way to get an Infinite Loch Ness monster from the Euclidean plane is using the operation well-known as the gluing straight segments. Actually, it consists of drawing two disjoint straight segments $l$ and $l^{'}$ of the same lengths on the Euclidean plane $\mathbb{E}^2$, then we cut along to $l$ and $l^{'}$ turns $\mathbb{E}^2$ into a surface with a boundary consisting of four straight segments (see Figure \ref{gluemarks}).

\begin{figure}[ht!]
 \begin{tabular}{cc}
    \includegraphics[scale=0.5]{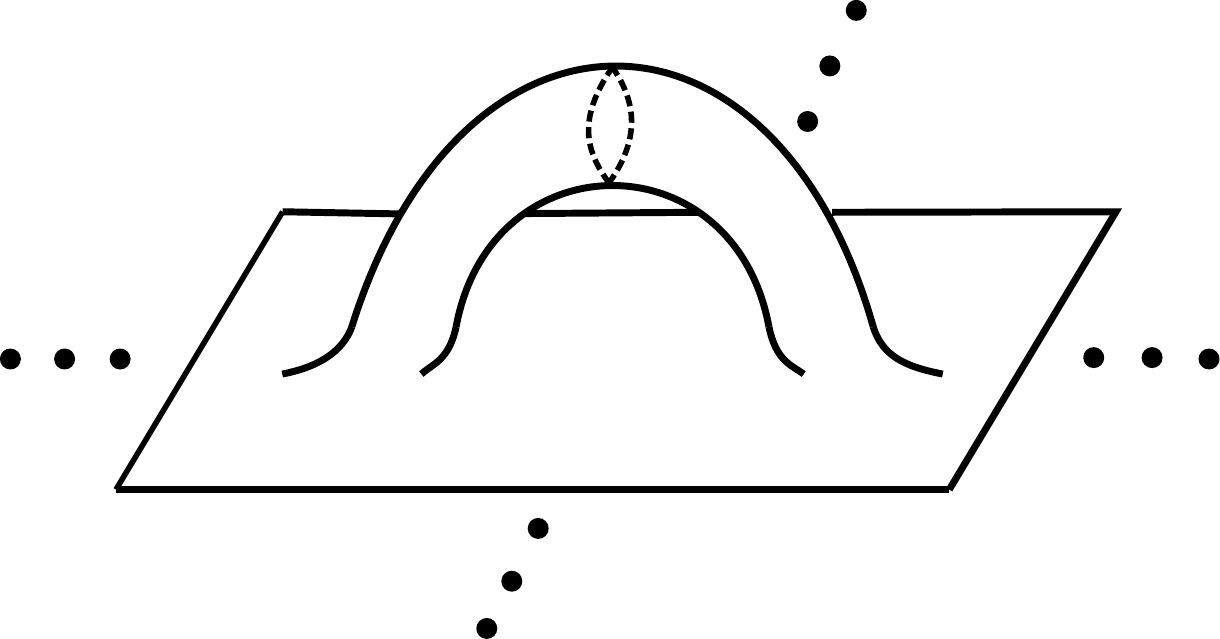}&    \includegraphics[scale=0.38]{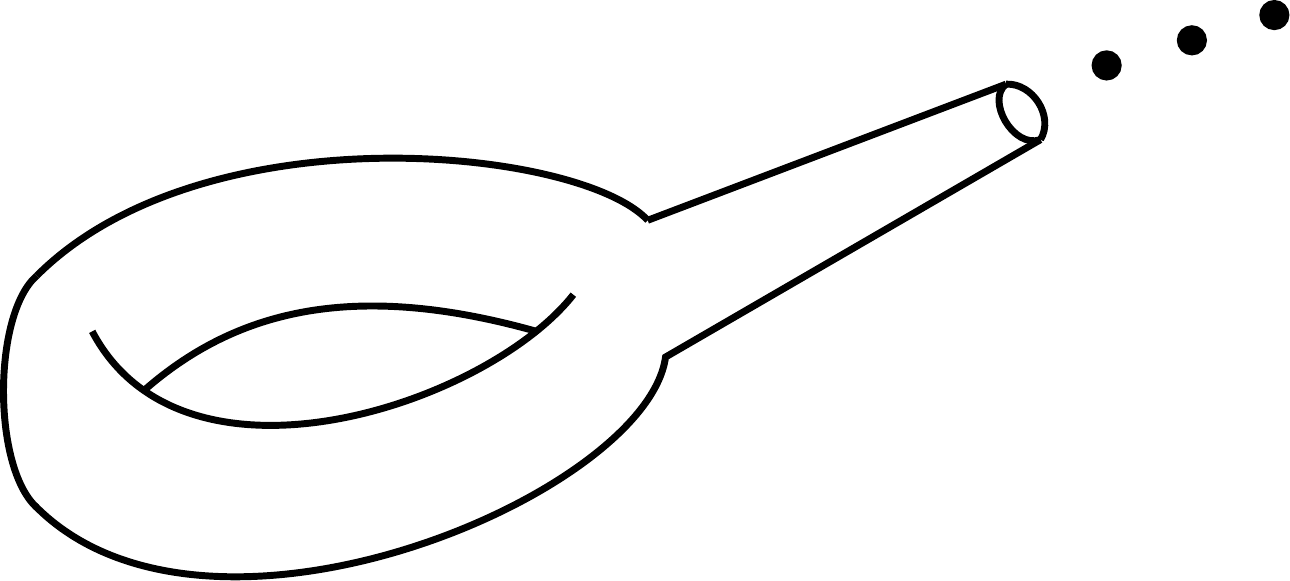} \\
   \emph{Gluing the two straight segments on $\mathbb{E}^2$.} &    \emph{Torus pictured by only one point.} \\
 \end{tabular}
 \caption{\textit{Gluing straight segments.}}
 \label{genus}
 \end{figure}

Finally, we glue this segments using translations to obtain a new surface $S,$ which is homeomorphic to the torus pictured by only one point (see Figure \ref{genus}). The operation described above is called \emph{gluing the straight segments} $l$ and $l^{'}$ \cite{RaVa}.

\begin{figure}[ht!]
  \centering
  \includegraphics[scale=0.5]{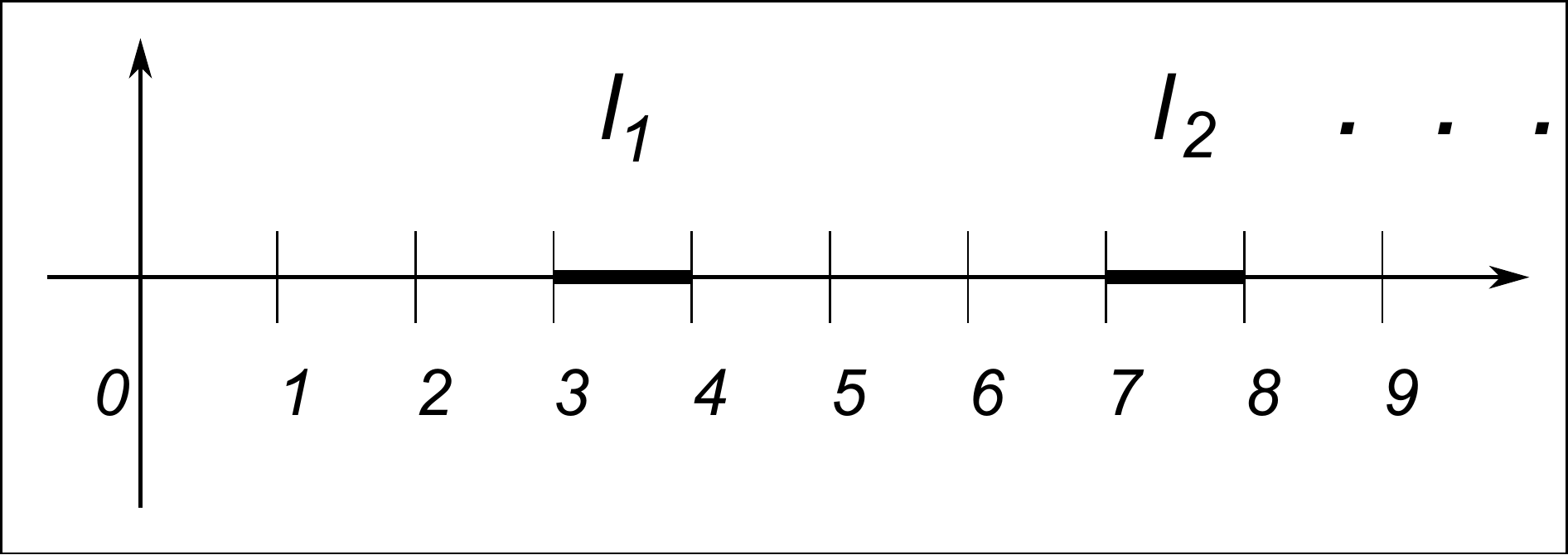}
  \caption{\emph{Countable family of straight segments $\mathcal{L}$.}}
  \label{glue}
\end{figure}

Note that to build a Loch Ness monster from the Euclidian plane using the gluing straight segments is necessary to draw on it a countable family of straight segment and suitable glue them.  It means, we consider $\mathbb{E}^2$ a copy of the Euclidean plane equipped with a fixed origin $\overline{0}$ and an orthogonal basis $\beta= \{e_{1},e_{2}\}$. On $\mathbb{E}^2$ we draw\footnote{Straight segments are given by their ends points.} the countable family of straight segments following:
\[
\mathcal{L}:=  \{l_{i} =((4i-1)e_{1}, \, 4ie_{1}) : \forall i \in \mathbb{N}\} \text{ (see Figure \ref{glue})}.
\]

Now, we cut $\mathbb{E}^2$ along $l_{i}$, for each $i\in \mathbb{N}$, which turns $\mathbb{E}^2$ into a surface with boundary consisting of infinite straight segments. Then, we glue the straight segments $l_{2i-1}$ and $l_{2i}$ as above (see Figure \ref{Figure3}). Hence, the surface $S$ comes from the Euclidean plane attached to an infinitely many handles, which appear gluing the countable disjoint straight segments belonged to the family $\mathcal{L}$. In other words, the mathematical object $S$ is the Infinite Loch Ness monster.

\begin{figure}[h!]
  \centering
  \includegraphics[scale=0.4]{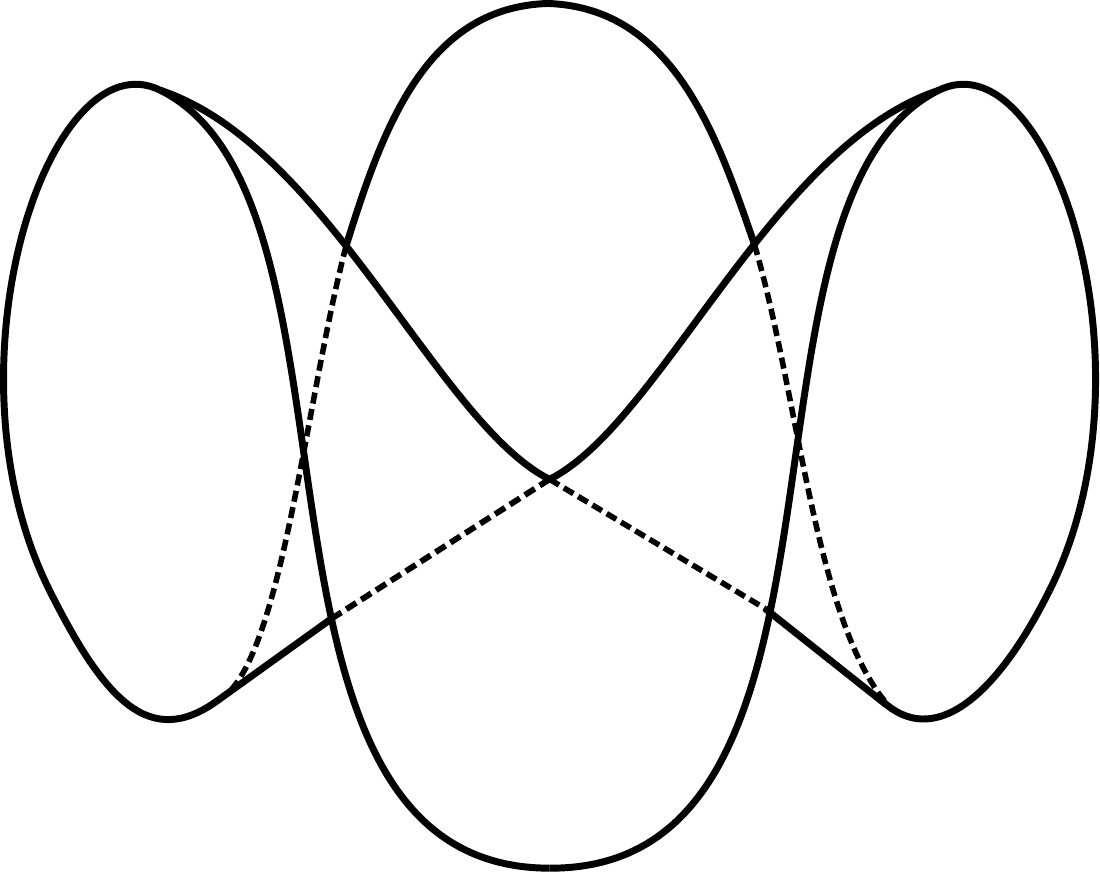}
  \caption{\emph{Cone angle singularity of angle $4\pi$.}}
  \label{cone}
\end{figure}

From view of differential geometric, the surface $S$ is conformed by two kind of points. The set of \emph{flat points} conformed by all points in $S$ except the ends of the straight segments $l_i$, for every $i\in \mathbb{N}$. To each one of this elements there exist an open isometric to some neighborhood of the Euclidean plane. Since the curvature is invariant under isometries then the curvature in the flat points is equal to zero.  The other ones, are called \emph{singular points}, in this case they are the end points of the straight segments $l_i$, for each $i\in\mathbb{N}$. Their respective neighborhood is isometric to cyclic branched covering $2:1$ of the disk in the the Euclidean plane, \emph{i.e.}, they are \emph{cone angle singularity of angle} $4\pi$ (see Figure \ref{cone}). The surfaces having this kind of structure are known as \emph{tame translation surfaces} (see \emph{e.g.}, \cite{PSV}).


\subsection{Hyperbolic Infinite Loch Ness Monster}
\label{hyperbolic}

An application of the Uniformization Theorem (see \emph{e.g.}, \cite{Abi}, \cite{Muc}) ensures the existence of a subgroup $\Gamma$ of the isometries group of the hyperbolic plane $Isom(\mathbb{H})$ acting on the hyperbolic plane $\mathbb{H}$ performing the quotient space $\mathbb{H} / \Gamma$ in a hyperbolic surface homeomorphic to the Infinite Loch Ness monster.  In other words, there exist a hyperbolic polygon $P\subseteq \mathbb{H}$, which is suitable identifying its sides by hyperbolic isometries to get the Infinite Loch Ness monster. An easy way to define the polygon $P$ is as follows\footnote{The reader can also found in \cite{AR} a great variety of hyperbolic polygons that perform hyperbolic surfaces having infinite genus.}.

\begin{figure}[h!]
  \centering
  \includegraphics[scale=0.4]{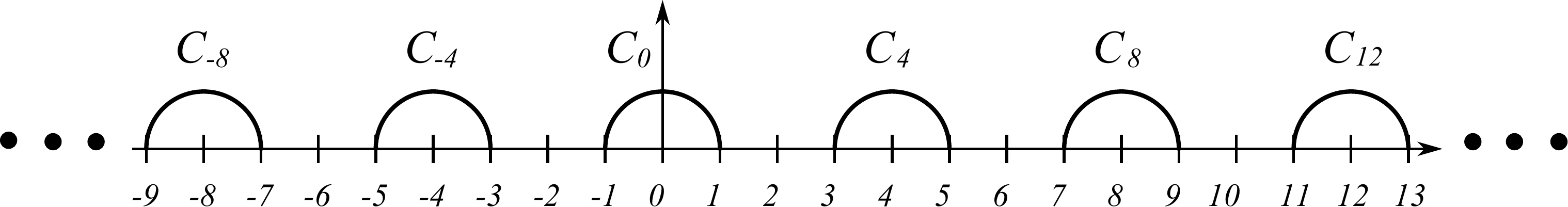}
  \caption{\emph{Family of half-circles $\mathcal{C}$.}}
  \label{circles}
\end{figure}

\begin{figure}[h!]
  \centering
  \includegraphics[scale=0.4]{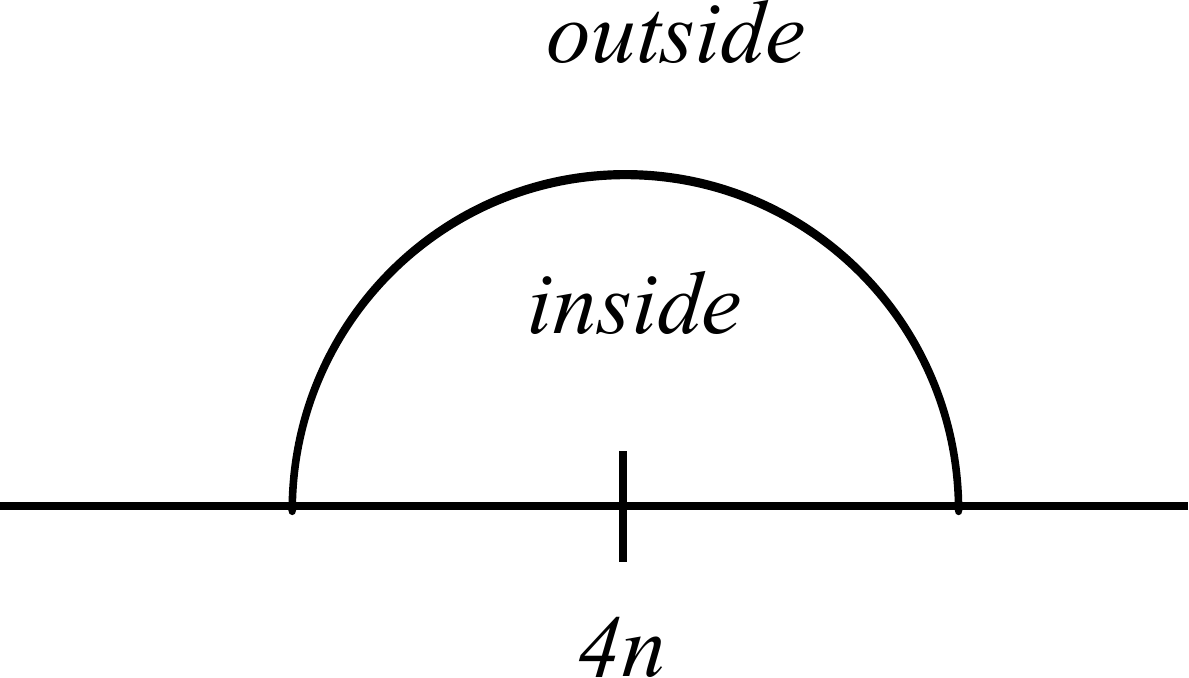}
  \caption{\emph{Inside and outside.}}
  \label{inside}
\end{figure}

\begin{figure}[h!]
  \centering
  \includegraphics[scale=0.4]{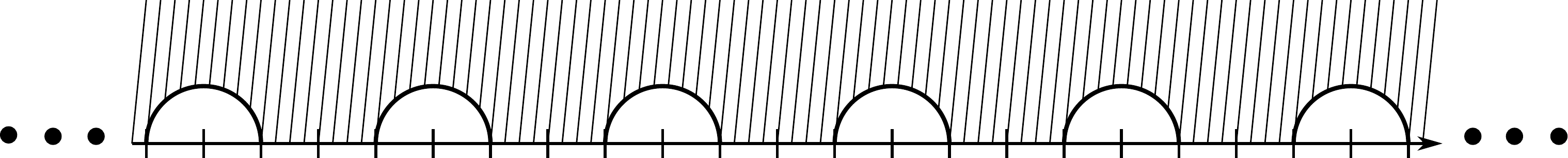}
  \caption{\emph{Family of half-circles $\mathcal{C}$ and hyperbolic polygon $P$.}}
  \label{poligono}
\end{figure}

First, we consider  the countable family conformed by the disjoint half-circles $\mathcal{C}=\{C_{4n}: n\in\mathbb{Z}\}$ having $C_{4n}$ center in $4n$ and radius equal to one, for every $n\in \mathbb{Z}$. See Figure \ref{circles}. In other words, $C_{4n}:=\{z\in \mathbb{H}: |z- 4n|=1\}$
 Removing the half-circle $C_{4n}$ of the hyperbolic plane $\mathbb{H}$ we get two connected component, which are called the \emph{inside} of $C_{4n}$ and the \emph{outside} of $C_{4n}$, respectively (see Figure\ref{inside}). They are denoted as $\check{C}_{4n}$ and $\hat{C}_{4n}$, respectively.
Hence, our connected hyperbolic polygon $P\subseteq\mathbb{H}$ is the closure of the intersection of the outsides following (see Figure \ref{poligono})

\begin{equation}\label{eq:4}
P:=\overline{\bigcap\limits_{n\in\mathbb{Z}} \hat{C}_{4n}}=\bigcap\limits_{n\in\mathbb{Z}} \{z\in\mathbb{H}: |z-4n|\geq 1\}.
\end{equation}

\begin{figure}[h!]
  \centering
  \includegraphics[scale=0.5]{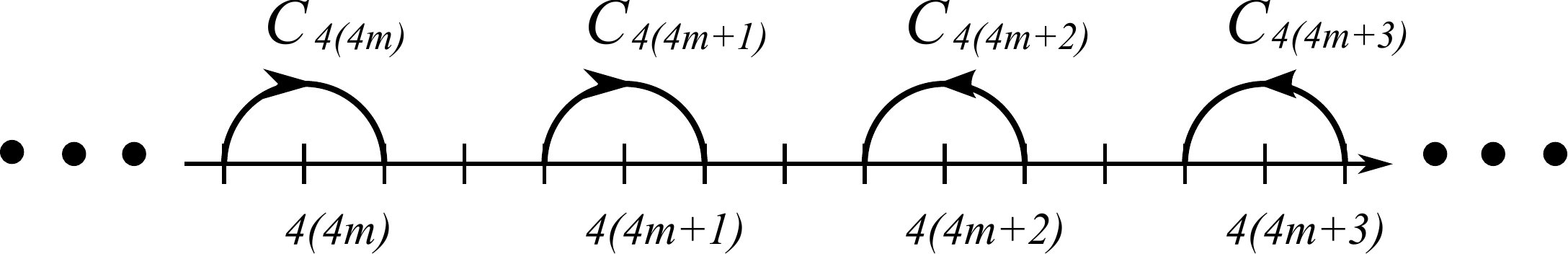}
  \caption{\emph{Gluing the side of the hyperbolic polygon $P$.}}
  \label{identificacion}
\end{figure}

The boundary of $P$ is conformed by the half-circle belonged to the family $\mathcal{C}$. Then for every $m\in\mathbb{Z}$ the hyperbolic geodesics $C_{4(4m)}$ and $C_{4(4m+2)}$ are identified as it is shown in Figure \ref{identificacion} by some of the following M\"obius transformations:

\begin{equation}\label{eq:5}
\begin{array}{rl}
f_{m}(z)        & :=\dfrac{(16m+8)z-(1+16m(16m+8))}{z-16m}\\
f_{m}^{-1}(z) & := \dfrac{-16mz+(1+16m(16m+8))}{-z+(16m+8)}.
\end{array}
\end{equation}

Analogously, the hyperbolic geodesics $C_{4(4m+1)}$ and $C_{4(4m+3)}$ are identified as it is shown in Figure \ref{identificacion} by the  M\"obius transformations:
\begin{equation}\label{eq:6}
\begin{array}{rl}
g_{m}(z) & :=\dfrac{(16m+8)z-(1+(16m+4)(16m+8))}{z-(16m+4)}, \\
g_{m}^{-1}(z) & := \dfrac{-(16m+4)z+(1+(16m+4)(16m+8))}{-z+(16m+8)}.
\end{array}
\end{equation}

\begin{remark}
Through the M\"obius transformations above, the inside of the half-circle $C_{4(4m)}$ (the half-circle $C_{4(4m+1)}$, respectively) is send by the map $f_m(z)$ (the map $g_m(z)$, respectively) into the outside of the half-circle $C_{4(4m+2)}$ (the half-circle $C_{4(4m+3)}$, respectively). Furthermore, the outside of the half-circle $C_{4(4m)}$ (the half-circle $C_{4(4m+1)}$, respectively)  is send by $f_{m}(z)$ (the map $g_m(z)$, respectively)  into the inside of the half-circle $C_{4(4m+2)}$ (the half-circle $C_{4(4m+3)}$, respectively).
\end{remark}

Hence, the hyperbolic surface $S$ get glued the side of the polygon $P$ is the Infinite Loch Ness Monster. From the polygon $P$ we deduce that noncompact quotient space $S$ comes whit a hyperbolic structure having infinite area. Fortunately, the identification defined above takes the pairwise disjoint straight segment in the boundary of $P$ performing into the only one end of the surface $S$.

\begin{figure}[h!]
  \centering
  \includegraphics[scale=0.4]{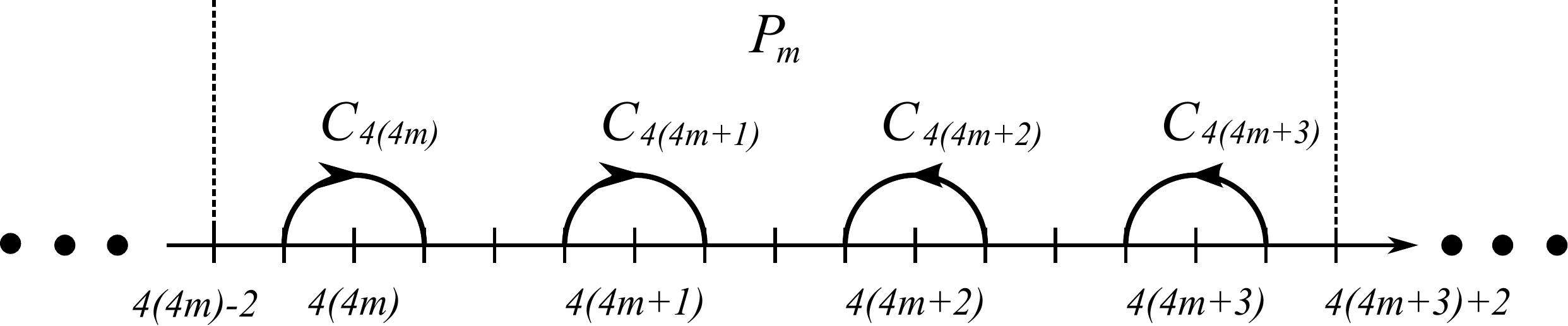}
  \caption{\emph{Subregion $P_m$.}}
  \label{subregion}
\end{figure}

\begin{figure}[h!]
  \centering
  \includegraphics[scale=0.4]{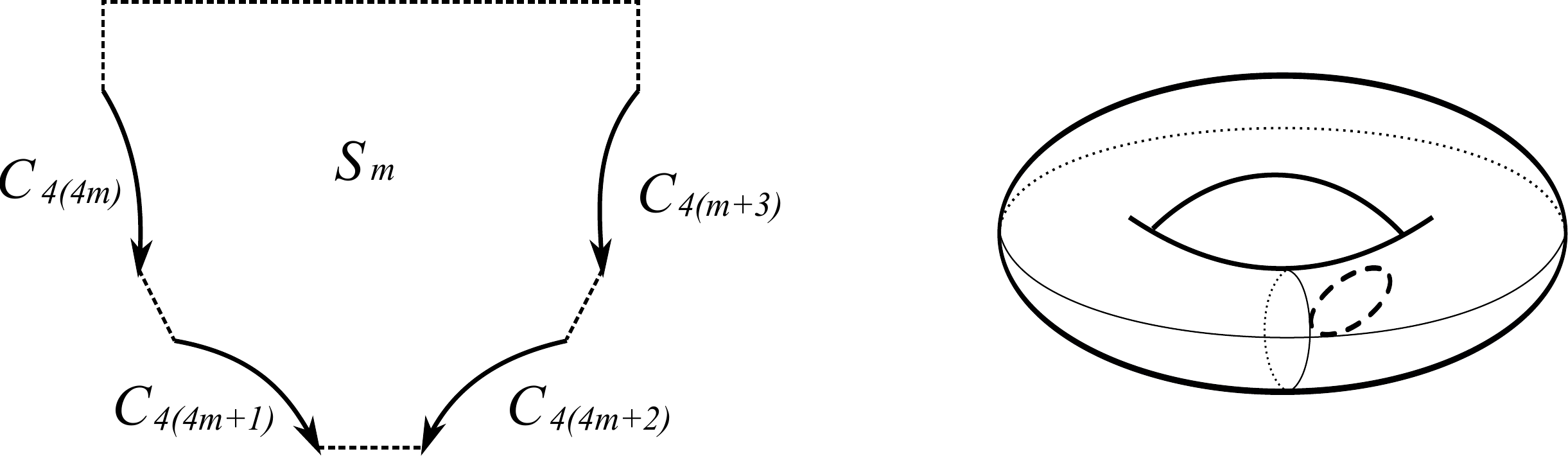}
  \caption{\emph{Topological subregion $P_m$ and torus with one hole $S_m$.}}
  \label{subsurface}
\end{figure}

Furthermore, for each integer number $n\in\mathbb{Z}$ we consider the subregion $P_m\subseteq P$, which is gotten by the intersection of $P$ and the strip $\{z\in \mathbb{H}:4(4m)-2 <Re(z)< 4(4m+3)+2\}$ (See Figure \ref{subregion}), then restricting to $P_m$ the identification defined above it is turned into a torus with one hole $S_m$ (see Figure \ref{subsurface}), which is a subsurface of $S$. Given the elements of the countable family $\{S_m:m\in\mathbb{Z}\}$ are pair disjoint subsurfaces of $S$ then it performs infinite genus in the hyperbolic surface $S$. In other words, $S$ is the Infinite Loch Ness monster.

From the analytic point of view, we have built a Fuchsian subgroup $\Gamma$ of $PSL(2,\mathbb{Z})$, where $\Gamma$ is infinitely generated by the set of M\"obius transformations $\{f_m(z), g_m(z), f^{-1}_m(z), g^{-1}_m(z): \text{ for all } m\in \mathbb{Z}\}$ (see equations \ref{eq:5} and \ref{eq:6}), having the subset $P\subseteq \mathbb{H}$ as fundamental domain\footnote{To deepen in these topics we suggest to reader \cite{MB}, \cite{KS}.  }. Then $\Gamma$ acts on the hyperbolic plane $\mathbb{H}$. Defining the subset $K\subseteq \mathbb{H}$ as follows,
\begin{equation}\label{eq:7}
K:=\{w\in\mathbb{H}: f(w)=w \text{ for any } f\in \Gamma-\{Id\}\}\subseteq\mathbb{H},
\end{equation}
the Fuchsian group $\Gamma$ acts freely and properly discontinuously on the open subset $\mathbb{H}-K$. Hence, the quotient space
\begin{equation}\label{eq:8}
S:= (\mathbb{H}-K)/\Gamma
\end{equation}
is a well-defined hyperbolic surface homeomorphic to the Infinite Loch Ness monster. Moreover, it follows from an application of the Uniformization Theorem that the fundamental group $\pi_1(S)$ of the Infinite Loch Ness monster is isomorphic to $\Gamma$.





\end{document}